\documentclass[11pt]{article}
\usepackage{amsmath,amsfonts,amssymb,amsthm,thmtools,thm-restate,mathabx,stmaryrd,enumerate,ytableau,hyperref,multirow}
\usepackage[boxed,linesnumbered]{algorithm2e}
\usepackage{cleveref}
\usepackage{fullpage}
\usepackage{color}
\usepackage{authblk}
\usepackage{float}
\usepackage{tikz}
\usepackage{ytableau}
\newif\ifdraft
\drafttrue
\DeclareFontFamily{OMX}{MnSymbolE}{}
\DeclareFontShape{OMX}{MnSymbolE}{m}{n}{
	<-6>  MnSymbolE5
	<6-7>  MnSymbolE6
	<7-8>  MnSymbolE7
	<8-9>  MnSymbolE8
	<9-10> MnSymbolE9
	<10-12> MnSymbolE10
	<12->   MnSymbolE12}{}
\DeclareSymbolFont{mnlargesymbols}{OMX}{MnSymbolE}{m}{n}
\SetSymbolFont{mnlargesymbols}{bold}{OMX}{MnSymbolE}{b}{n}
\DeclareMathDelimiter{\llangle}{\mathopen}{mnlargesymbols}{'164}{mnlargesymbols}{'164}
\DeclareMathDelimiter{\rrangle}{\mathclose}{mnlargesymbols}{'171}{mnlargesymbols}{'171}

\makeatletter
\@addtoreset{equation}{section}
\makeatother

\title{Tri-skill variant Simplex and strongly polynomial-time algorithm for linear programming}
\author{P. Z. Wang\textsuperscript{\rm 1}\thanks{Liaoning Technical University, Fuxin, Liaoning, 123000, China. Email:peizhuangw@126.com}, J. He\textsuperscript{\rm 2},  H. C. Lui\textsuperscript{\rm 1}, Q. W. Kong\textsuperscript{\rm 3}, Y. Shi\textsuperscript{\rm 4}, S. Z. Guo\textsuperscript}
\affil[1]{Institute of Intelligence Engineering and Math, Liaoning Technical University, Fuxin, Liaoning, China}
\affil[2]{School of Software and Electrical Engineering, Swinburne University of Technology, Hawthorn, Australia}
\affil[3]{Institute of Information Engineering, Nanjing University of Finance and Economics, Nanjing, China}
\affil[4]{Key Laboratory of Big Data. Chinese Academy of Sciences, Beijing, China}
\date{}
\begin{document}
	\maketitle
	
	\begin{abstract}
		The existence of strongly polynomial-time algorithm for linear programming is a cross century international mathematical problem, whose breakthrough will solve a major theoretical crisis for the development of artificial intelligence. In order to make it happen, this paper proposes three solving techniques based on the cone-cutting theory: 1. The selection of cutter: principles highest vs. deepest; 2. The algorithm of column elimination, which is more convenient and effective than the Ye-column elimination theorem; 3. A step-down algorithm for a feasible point horizontally shifts to the center and then falls down to the bottom of the dual feasible region $D$. There will be a nice work combining three techniques, the tri-skill is variant Simplex algorithm to be expected to help readers building the strong polynomial algorithms. Besides, a variable weight optimization method is proposed in the paper, which opens a new window to bring the linear programming into uncomplicated calculation.\\
		
		$\mathit{Keywords}:$ Linear programming, Strongly polynomial-time algorithm, Cone-cutting,
		Column elimination algorithm, Factor space.
	\end{abstract}
	
	\section{Introduction} \label{sec:introduction}
	Linear programming is not only an indispensable computing tool for intelligent decision making and data science, but also a potential cornerstone of intelligent science with its causal interaction and dialectical connotation\cite{1}. It is also an important subject concerned by the application of factor space theory\cite{2,3}. The Simplex algorithm proposed by G. B. Dantzig\cite{4} is a precious pearl in mathematics, in mathematics, which is not only beautiful in mathematical theory, but also has been widely used in practice. Although Klee-Minty counterexample\cite{5} points out that the deepest Simplex may appear exponential explosion, but people always hope to find generalized Simplex polynomial algorithm. Owing to distinguish it from the Karmarkar's weak polynomial algorithm\cite{6}. The researchers are still eager to find the strongly polynomial algorithm still. The breakthrough of the cross-century mathematical problem\cite{7}, will save the development of artificial intelligence from an important theoretical crisis. In 1989, Y. Ye proposed the elimination theorem\cite{8}, which is of great significance. Under the influence of Z. Z. Zhang\cite{9}, P. Z. Wang proposed the cone cutting theory\cite{10,11} in 2014, which provided a clear geometric description for the Simplex algorithm in the dual space. On this basis, we proposed the gravity sliding algorithm\cite{12,13}. Based on the mentioned works, the work of this paper is to continue to put forward three solving techniques for linear programming under the guidance of this geometric vision. Firstly, the highest principle for the selection of cutting plane from the special t-value tableau such that the elevation of the new cone vertex is changed as higher as possible. The second is the new column elimination rule, which is more easy to use and more effective than Ye's column elimination theorem. The third is making a feasible point to horizontally shifts to the feasible center and falls down to the bottom of the feasible region $D$. These three techniques all have the function of accelerating solution, and can be matched with each other to produce joint effect, the tri-skill is variant Simplex algorithm to be expected to  help readers building the strong polynomial algorithms. 
	
	Structure of this paper is as follows: The second section introduces the theory of cone cutting, which is the theoretical basis of this paper. It is quoted from the paper \cite{13}, but it is updated. The third section introduces the highest principle for the selection of cutting plane, which is a discretization of the gradient method in the non smooth convex domain. Using this principle in the Klee-Minty counter
	example, the solution can be obtained in one step. Section 4 introduces a new column elimination
	theorem. The condition of column elimination is that the norm vector of the cutting plane must be non-negative or non-positive (homo-symbolic). We provides a homo-symbolic theorem to convert a column
	vector mixed with positive and negative coefficients into a homo-symbolic vector, and into the new column elimination theorem. The fifth section introduces how to find a feasible interval on a ray, and conveniently solves the difficulty of horizontal edge: if $0$ appears in the right column of the Simplex table, then the base row is fixed, the Simplex method will not work well here. In fact, as long as we find out a feasible point on this edge, it must be an optimal point, and then the problem is solved. If there is no point, the edge must be completely cut off. Instead of avoiding the meeting of horizontal edge, we should actively capture them. The sixth section introduces the third technique, which constructs a horizontal plane passing through a given feasible point and extents the feasible rays towards to each edge and forms an m-vertices polygon. Then horizontally shifts the feasible point to the center of the polygon, the center of this polygon, and performs gravitational descend to the bottom of the feasible region.  The deeper the falling, the faster the solbing.. The seventh section emphasizes that these three techniques can be combined together to create synergistic effect. The readers can make use of these techniques to develop the strongly polynomial time algorithm. Finally, a short conclusion is given in eighth section.
	\section{Cone-cutting theory} \label{sec:Cone-cutting theory}
	\subsection{Mathematical definition of cones}\label{sec:Mathematical definition of cones}
	A non-empty set $F$ in the space $Y=R^{m}$ is called an affine set, if for any two points $P, Q\in F$ and $t\in (-\infty, +\infty)$, there is $P+t(Q-P)\in F$. Two affine sets $F$ and $F'$ are called parallel, if there is $P \in F$  such that $P+Q\in F'$ is true for any $Q\in F'$, and there is $Q\in F$ 'such that $P+Q\in F'$ is true for any $P\in F$. An affine set is a subspace if and only if it contains the origin $O$. Any affine set is parallel to and only parallel to a subspace, so the dimension of affine sets can be determined. A $(m-1)-$ dimensional affine set in $Y$ is called a hyperplane, hereafter referred to as plane.
	
	Linear inequality constraints in linear programming problems can be transformed into equality constraints in terms of boundary surface, but the direction of inequality must be defined: In the search that maximizes the target value, only ``$\leq $'' is allowed, while in the search that minimizes the target value, only ``$\geq $'' is allowed. In minimum optimization, if the inequality $y_{1}-y_{2}\leq 1$ occurs, we must rewrite it as $-y_{1}+y_{2}\geq -1$, and then write the equation $-y_{1}+y_{2} = -1$. The norm vector to this equation must be $(-1, 1)$ and cannot be $(1, -1)$. The norm vector of all the planes is only free to be multiplied by a positive number and not by a negative number, so that the plane is going to be directed. A plane is determined by its norm vector $\tau$ and constant $c$ for simple, the plane is denoted as $(\tau_{j})$.
	
		Let $y=(y_1,...,y_m)$ be a row vector and $\tau=(\tau_1,...,\tau_m)^T$ be a column vector, then a plane is an $m-1$ dimensional affine set constrained by an equation $y\tau =c$, where $c$ is a constant and $\tau$ is referred as the normal vector. Furthermore, a facet denotes an inequality constraint $y\tau \geq c$. For simplicity, we use the symbol $(\tau)$ to denote both a plane and a facet. A point $P$ in $Y$ is called accepted by the facet if $y\tau \geq c$. We also say that it lies inside the interior of the facet. Otherwise, $P$ is rejected by the facet as it lies on the exterior of the facet.
	
	There are many kinds of cones, which are discussed in details in the combinational optimization textbook\cite{14}, but the cone we say here is only the special one. In order to get rid of complexity, this paper gives the following definitions of cone directly and concisely instead of using deeper terminologies.\\
	\textbf{Definition 2.1} \cite{11}\label{definition2.1}
	Let $y\tau_{j}=c_{j} (j=1,...,m)$ be a group of planes, the norm vector of them form the matrix $B=(\tau_1,...,\tau_m)$. If the rank $r(B) = m$, then their common accepted area is called a cone, denoted as:
		\begin{center} $C=\{ y \in Y | y\tau_1\geq c_1,...,y\tau_1\geq c_m\}$\\ \end{center}
		
		$B=(\tau_1,...,\tau_m)$ is called the face matrix of $C$ and each plane $y\tau_{j}=c_j$ is called a face of $C$.\\
		Note that $y$ is a row vector, $\tau_{j}$ is a column, therefore $y\tau_{j}$ is the inner product $(y^T, \tau_{j})$.\\
		Since $B$ is full rank, there must be a unique intersection point $V=V_{1\times m}$, called the vertex of cone $C$. It is obvious that:
		\begin{equation}\label{2.1}
			V=c_{B} B^{-1}
		\end{equation}
		
		Where $c_{B}=(c_1,..., c_m)$.
		
		In the $m$-dimensional space, any $m-1$ faces must yield a straight line as long as the rank of the matrix formed by their norm vectors is $m-1$. Let $L_i$ be the intersection line of all faces except the $i$-th face, which is called the $i$-th edge-line of $C$. It is easy to prove that if $r(B)=m$, then $L_1,…,L_m$ can not lie in a common face, so that the $i$-th edge-line can not lie in the $i$-th face. All edge-lines intersect in the vertex $V$. The ray being the half of $L_i$ within the cone $C$ is called the $i$-th edge of $C$, denoted as $\underline{L_i}$. Its direction is denoted as $e_i$, a row vector. Since a direction has no fixed length, set that:
		\begin{equation}\label{2.2}
		e_i^{*}=e_i /e_i\tau_{i}=e_i / (e^{T}_i,\tau_{i}) \; (i = 1,...,m),  
		\end{equation}
		
		$E^{*}=(e^{*}_{1},...,e^{*}_{m})^{T}$ is called the regular edge matrix. $E^{*}$ can be determined directly by the basis matrix $B$. Obviously we have that:
		\begin{equation}\label{2.3}
			e^{*}_i\tau_{i} = 1 \; (i= 1,...,m).
		\end{equation}
		The position of face and edge is equal and they can be determined mutually. Since the $i$-th edge is the intersection of all except the $i$-th face, the $i$-th face should also be spanned by all except the $i$-th edge. We can define a cone in terms of edges:\\
		\textbf{Definition $2.1^{\prime}$} \cite{11}\label{defi2.1'}
		Given a point $V$ and matrix $E=(e_1,...,e_m)^{T}$. If $r(E)=m$, denote that
		\begin{center}
				$C=$\{y $\in$Y | $\exists$($\lambda_{1},...,\lambda_{m})\geq$ 0, $\lambda_{1}+$,...,+$\lambda_{m}=1$, such that $\lambda_{1}e_{1}+,...,+\lambda_{m}e_{m}=y$\}.
		\end{center}
	
		Which is called an $m$-cone, $V$ is called the vertex of the cone, $\underline{L_i}=V+t e_i (t\geq 0)$ is called the $i$-th edge of $C$.
		
		We are going to prove the equivalence between definition \ref{definition2.1} and definition \ref{defi2.1'}:\\
		\textbf{Theorem 2.1} \label{theo2.1}The regular edge matrix $E^*$ and face matrix $B$ are reciprocal: $E^*=B^{-1}$.\\
		\textbf{Proof} For any $i =1,...,m$, $L_i$ is the intersection line of all faces except the $i$-th one. Since $L_i$ is on all these faces, it is perpendicular to the norm vector of all these faces:
		\begin{equation}\label{2.4}
			e^*_i \tau_{j} = 0 \; (j\neq i). 
		\end{equation}
		
		According to \ref{2.3} and \ref{2.4}, we have that $E^*B=I$. \; \; \textbf{End}\\
		\textbf{Definition 2.2} If the face matrix of $C_2$ is the regular edge matrix of $C_1$, then, cone $C_2$ is called the anti-cone of cone $C_1$, denoted as $C_2 = C_1^{-1}.$
		
		Anticonicity is obviously symmetric:\\
		\textbf{Corollary 2.1} If cone $C_2$ is the anticone of cone $C_1$, then cone $C_1$ must also be the anticone of cone $C_2$.\\
		\textbf{Proof} Assume that the face matrices of $C_1$ and $C_1$ are $B_1$ and $B_2$, and the edge matrices of them are $E^*_1$ and $E^*_2$ respectively. Since cone $C_2$ is the anticone of cone $C_1$, $B_2= E^*_1$. Take the inverse of both ends of this equation, and it can be known from Theorem \ref{theo2.1} that $E^*_2 = B_1$, the face matrix of $C_1$ is the regular edge matrix of $C_2$. Therefore, $C_1$ is anticone of cone $C_2$. \textbf{End}
		
		It is obvious that:
		\begin{equation}\label{2.5}
				(C^{-1})^{-1}=C.
		\end{equation}
		\subsection{Cone interpretation of Simplex table}\label{sec:Cone interpretation of Simplex table}
		Given a problem of linear programming:
		\begin{center}
			$(P)\; max\{cx|Ax\leq b\}$,
		\end{center}
	
		Where $x=x_{n\times 1}$ and $b=b_{m\times 1}$ are column vectors, $y=y_{1\times m}$ and $c=c_{1\times n}$ are row vectors, $A=A_{m\times n}$. The original Simplex tableau is given as follows:
		\begin{center}
		Table 2.1 Original Simplex table
		\begin{equation}\label{Tab2.1}
		\bordermatrix{%
		&x_1 & \cdots &  x_n & y_1 & \cdots & y_m &\cr
		&\alpha_{11} & \cdots & \alpha_{1n} & 1 & \cdots &0& b_{1}\cr
		&\vdots & \ddots & \vdots & \vdots &\ddots &\vdots &\vdots\cr
	    &\alpha_{m1} &\cdots& \alpha_{mn} &0& \cdots&1& b_{m}\cr
     	&c_{1} &\cdots& c_{n} &0& \cdots&0& 0\cr }  \tag{$T_o$}
	    \end{equation}
	\end{center}

		Denote $A^+=(A, I)$, the header $(x_1,..,x_n; y_1,...,y_m)$ is collectively denoted as $(\tau_1,..., \tau_{n+m}$, Let $c^{+}=(c_1,...,c_n,0...0)$. The Cone-Cutting theory interprets this original tableau as follows: The identity matrix $I$ represents an original cone $C_o$, Its vertex is the origin of the coordinate system $O=(0,...,0)$; its edges are the coordinate axes. Each row of $I$ is the directional coefficients of the edge. The row vector below this matrix $I$ is the coordinate of the vertex of this cone. We can combine matrix $I$ and the vertex vector below it together to form a matrix $\underline{I}$. This is called the initial cone matrix $C_o$. Let $\underline{A}$ denotes the extension of $A$ to include the row vector $c$ below it, then the equation $yA_j=c_j$ represents a dual constraint plane, called a cutting plane. All cutting plane cuts the $Y=R_m$ space into a feasible region $D$. The points inside $D$ are called the feasible points.  
		
	In this paper, we say that point $P_A$ is lower than point $P_B$ if $P^{T}_Ab<P^{T}_Bb$, where $b$ is the dual target vector.
	
	That is, all measurement comparison is relative to the $b$ vector; and we refer $u=P^{T}_Ab$ as the height of $P_A$.
		
		The unit square $I$ represents the original cone $C_o$, whose vertex is the origin $O=(0,...,0)$, the edges
		are the axes. The unit square $I$ is the edge matrix, each row is the direction of an edge, and the row
		vector below $I$ is the coordinate of the vertex of the cone. Let $\underline A$ be the combining $A$ with its lower row, then, $I$ is called the cone matrix of tableau \ref{Tab2.1}. In this paper, we always assume that the dual target vector is non-negative, that is, $b\geq 0$. Under the assumption that vertex $O$ is in its bottom of $C_o$, if cone $O$ is accepted by all constraint planes, it is the optimal point! Otherwise, there must be a constraint plane
		cut off $O$, and form a newer cone $C$ containing the feasible region $D$. The algorithm of cone cutting
		ensures that the new cone must be regular, $i.e.$ the vertex of cone is in its bottom. Repeat, once the cone
		vertex cannot be cut, it is an optimal point, this is the thought of the cone cutting.
		
		Set the index set of all constraint planes as $J=\{1,...,n+m\}$, let $\underline{B}$ be a subset of $J$ , containing m indices. If the indicated matrix $B$ is full rank, then $\underline{B}$ is called a base index set. The Simplex theory rely on the base transformation:
		\begin{equation}\label{2.6}
			y'=yB^{-1}+V, y=(y'-V)B.
		\end{equation}
		
		Suppose that the cutting plane ($\tau_{j}$) has equation $y'\tau_{j}=c_j$. According to \ref{2.6}, its equation can be written as $(y B^{-1}+V)\tau_{j}=c_{j}$, or $y B^{-1}\tau_{j}=c_j-V\tau_{j}.$ when expressed in the transformed coordinate system. Let $c_j^{\land}=c_j-V\tau_{j}$.\\
		\textbf{Definition 2.3} The constant $c_j^{\land}$ is called the cutting degree of plane ($\tau_{i}$) with respect to the vertex $V$. Vector $c^{\land}=(c^{\land}_1,...,c^{\land}_{n+m})$ is called the cutting vector.
		
		Cutting degree $c^{\land}_j$ measures the norm distance from a cutting plane to cone vertex; $V$ is cut off by ($\tau_{j}$), if and only if $c^{\land}_j>0$. cut degree is the transformation of the constant term in the equation of ($\tau_{j}$):
		\begin{equation}\label{2.7}
				c^{\land}_j=c^{+}_j-c_BB^{-1}A_j
		\end{equation}
	
		Where $c_B$ is the subset of $c^+=(c_1,..., c_m,0,...,0)$ with index set $\underline{B}$.\\
		\textbf{Proof}. For $c^{\land}_j =c^+_j-V \tau_{j}$, \ref{2.7} is true by means of \ref{2.1} . \textbf{End}\\
		\textbf{Corollary} When $c^{+}_j=0$, We have that:
		\begin{equation}\label{2.8}
			c^{\land}_j=-V\tau_{j}.
		\end{equation}

		There are many papers use geometry to describe LP solution, but descriptions could not clearly
		stated in space $X$; Cone cutting theory describes pivoting in the dual space $Y=R^m$ for the dual problem:
		\begin{center}
			$(D) \; min \; \{yb|yA \geq c, y \geq 0\}$
		\end{center}
	
There are $n+m$ constraint planes, the new usage of symbol $\tau_j$ now stands for the norm vector of the plane $(\tau_{j})$. Its initial vector is that $\tau_j^{o}j=A_j^{+} (j=1,...,n+m)$. After the basis transformation, the general Simplex
		matrix $T=(\tau_{ij})$ is written in the following tableau:
		\begin{center}
			Table 2.2 Transformed Simplex tableau
			\begin{equation}\label{Tab2.2}
			\begin{pmatrix}
			B^{-1}A	& B^{-1}  & B^{-1}b \\ 
			c-c_BB^{-1}A	& -c_BB^{-1} & -c_BB^{-1}b
			\end{pmatrix} \tag{$T$}
			\end{equation}
		\end{center}
	
		When $\underline{B}=\{n+1,...,n+m\}$, $B$ returns to the original base $I: B=I$, then $c_B=0$, so the tableau \ref{Tab2.2} returns to \ref{Tab2.1}.\\
		\textbf{Theorem 2.1} enables us to replace $B^{-1}$ with $E^*$ in Simplex tableau \ref{Tab2.2}, thus obtain the cone interpretation of Simplex tableau:
		\begin{center}
			Table 2.3 Cone-Simplex tableau
			\begin{equation}\label{Tab2.3}
			\begin{pmatrix}
			E^{*}A	& E^{*}  & s=E^{*}b \\ 
			c-VA	& -V & -h(V)
			\end{pmatrix} \tag{$T$}
			\end{equation}
		\end{center}
	
		Where $s=E^{*}b$ is called the slop vector of edges: $s_i=0$ means that the edge $\underline{L_i}$ is horizontal with
		respect to $b$; $\underline{L_i}$ going up (down) when $s_i > 0 \;(s_i <0)$.
		
		We have to emphasize that the main matrix can be regarded as a couple:
		\begin{center}
		$\begin{pmatrix}
		E^{*}A \\
	c^{\land}=c-VA
		\end{pmatrix}\rightarrow A'$\;\;\;\;\;
	$	\begin{pmatrix}
		E^*\\
		-V
		\end{pmatrix}\rightarrow C$
	\end{center}

		The matrix $\underline{A'}$ is transferred from $A$ in table \ref{Tab2.1}, which is called the constraint matrix, the $j$-th column shows the new equation $y\tau_{j}=c^{\land}_j$ in stead of equation $yA_j=c_j$. Note that $A'=E^{*}A$. It means that for any $\tau_{ij}$ in table \ref{Tab2.3}, we have that
		\begin{equation}\label{2.9}
		\tau_{ij}=e_i^{*}A_j,
		\end{equation}
		
		For $B$, each row vector in $E^*$ is the edge vector $e^*_i$ of a cone and the last entry of this sub-matrix is the negative value of the vertex coordinates. Thus the submatrix $B$ describes the cone completely and is called the cone matrix.
		For $C$, $s$ is the edge slope vector as mentioned before. Each element $s_i=e^*_ib$ meansures the slope of the edge vector $e^*_i$ against the target vector $b$. The last entry is the negative value of $h(V)=Vb$. which is actually the current value of the objective function that we want to minimize.
		
		Thus, each term in the table has a clear geometric interpretation in cone cutting theory. This is
		the new vision of the Simplex method brought by cone cutting theory.
		\subsection{Cone-cutting vs Pivoting}
		Whenever a vertex $V$ of cone $C$ is cut by a plane ($\tau_{j^*}$). How do we generate a new cone $C'$? Where is its vertex $V'$, and what are directions of new edges? There is an algorithm Cone-cutting to resolve these
		problems. This algorithm is not described in detail in this paper. Because, the significance of cone cutting
		lies not in the very algorithm but in its geometric description for Simplex. In a nutshell, Let $Q_i=V+t_{ij} e_i$ be the intersection of plane ($\tau_{i}$) and edge $\underline{L_i}$, it must hold that
		\begin{equation}\label{2.10}
		t_{ij}=	c^{\land}_j/\tau_{ij}(\tau_{ij}\neq 0)
		\end{equation}
		
		$Q_i$ is called real, if $t_{ij}\geq 0$; it is called virtual, If $t_{ij} < 0$, a virtual $Q_{ij}$ is located on the inversed ray of $L_i$. The vertex of the new cone is that:
		\begin{equation}\label{2.11}
		V’=Q_{i^*} ; \; i^* = \mathop{\arg\min}_{i}\{t_{ij} e_i b| t_{ij}>0\}
		\end{equation}
		
		Which is the lowest real intersection of edges with the cutter ($\tau_{ij}$).\\
		The edges’ directions of $C’$ are calculated as follows:
		
		(1) $e^{'}_{i^*}=e_{i^*}$ ;
		
		(2) If $i\neq i^*$, and $t_{ij}>0,$ then $e^{'}_i=Q_i-Q_{i^*}$ ;
		
		(3) If $i\neq i^*$, and $t_{ij}<0,$ then $e^{'}_i=Q_{i^*}-Q_i$ ;
		
		(4) If edge line $L_i$ parallels to plane $\tau_{j^*}$, then the $e^{'}_i=e_i$.
		
		Cone cutting demonstrates the pivoting process, pivoting is to bring in a non-basic variable $x_{j^*}$ to replace a basic varible $x_{i^*}$. Once $j^*$ is determined, $i^*$ is computed according to:
		\begin{equation}\label{2.12}
		i^{*}=\mathop{\arg\min}_{i}\{s_i/\tau_{ij^*}|\tau_{ij^*} > 0, i=1,...,m \} 
		\end{equation}
		
		Simplex performs the pivoting operation at the pivoting point $(i^*, j^*)$. For Cone-Cutting, a non-basic facet $j^*$ not in $B$ is selected to replace a basic facet $i^*$ from the set of basic facets. Referring to \ref{2.12}, since $t_{ij}=c^{\land}_j/\tau_{ij}$, that $\tau_{ij}=c^{\land}_j/t_{ij}$, while $s_i=e_ib$, we have that $s_i/\tau_{ij}=t_{ij} e_i b/c^{\land}_j$. Note that $c^{\land}_j$ is a positive
		number, which has no relation with the index $i$, so that \ref{2.11} equivalents to \ref{2.12}. It means that cone cutting is the very pivoting.
		
		Adding a column to the leftest of the Simplex table to represent $\underline{B^T}$. The initial base indices was $\underline{B_o}=\{n+1,..., n+m\}$ and the base faces was $(y_1),...,(y_m)$, $\underline{B}$ changes a name during each pivoting; the name of expelled face can be find out from $\underline{B}$ at $i^*$-th row.\\
		\textbf{Example 2.1} Given a linear programming problem:
		\begin{center}
			$Max\; \; 2x_1+x_2$\\
		$\; s,t.\; x_1-x_2 \leq 1;$\\
		$\; \; x_1+x_2 \geq  2.$\\
		$\;\; -x_1 \leq 0,-x_2 \leq 0$
		\end{center}

		Its standard Simplex table \ref{Tab2.1} can be written as follow Table 2.4:
		\begin{center}
			Table 2.4 Simplex table of example 2.1
			\begin{equation}\label{Table 2.4}
			\bordermatrix{%
				& x_1       & x_2     & y_1    &y_2\cr 
				x_1    & 1         & -1       &1     & 0&1\cr 
				y_2    & 1        & 1       &0     & 1&2\cr 
				c^{\land} & 2    &1  &0     &0 &0\cr 
			}\notag
			\end{equation}
		\end{center}
	
\begin{figure}[htpb]
			\centering
			\includegraphics[width=0.25\linewidth]{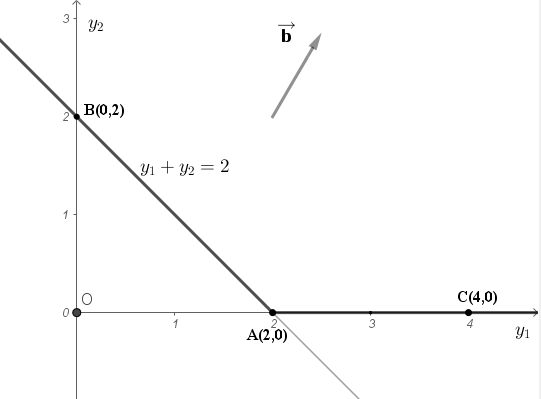}
			\caption{Cone $C_1$}
			\label{fig1}
		\end{figure}

		We can see from Figure \ref{fig1}, the original cone $C_o$ is the first quadrant of $Y=R^2$, the vertex of the cone is the
		origin $O$, and its two edges are the axis vectors $e_1=(1,0)$ and $e_2 =(0,1)$, which are written in the cone
		matrix $C$ of \ref{Tab2.1} respectively. The dual target vector is $b=(1,2)^T$ and $O$ is the lowest point of the first
		quadrant. Since the first number of $c^{\land}$ is a positive number 2, it means that face $y_1+y_2=2$ (the line $(x_1)$
		in the figure) can cut off the vertex $O$.
	
		Taking $j^*=1$, according to \ref{2.12}, we have $i^*=1$; Doing pivot at (1,1) get table as follows:
		\begin{center}
		\begin{tabular}{cccccc}
		$\underline{B}$	& $x_1$ &  $x_2$&  $y_1$& $y_2$ &  \\ 
		$x_1$	& 1 & -1 &1  &0  &1  \\ 
			$y_2$& 0 &2  &-1  &1  & 1 \\ 
		$c^{\land}$	& 0 & 3 &-2  &0  & -2 \\ 
		\end{tabular} 
		\end{center}
		
		In the left column of the table, cutter $(x_1)$ replaces face $(y_1)$. The vertex of the new cone is shown under the cone matrix: $V_1=-(-2,0)=(2,0)$, and the new cone is angle $\angle BAC$, where $B=(0, 2)$, $A=(2, 0)$ and $C=(4, 0)$. Two edge directions are listed in the second table: $e'_1=(1,0), e'_2=(-1, 1)$. Since $c^{\land}_2 = 2>0$, $V_1$ can be cut by $(x_2)$. Take $j^*=2$, according to \ref{2.12}, $i^*=2$, doing pivot at \ref{2.2} and get the following table:
			\begin{center}
			\begin{tabular}{cccccc}
				$\underline{B}$	& $x_1$ &  $x_2$&  $y_1$& $y_2$ &  \\ 
				$x_1$	& 1 & 0 &0.5  &0.5 &1.5 \\ 
				$y_2$& 0 &1  &-0.5  &0.5  & 0.5 \\ 
				$c^{\land}$	& 0 & 0 &-1.5  &-0.5  & -3.5 \\ 
			\end{tabular} 
		\end{center}
	
\begin{figure}[htpb]
			\centering
			\includegraphics[width=0.25\linewidth]{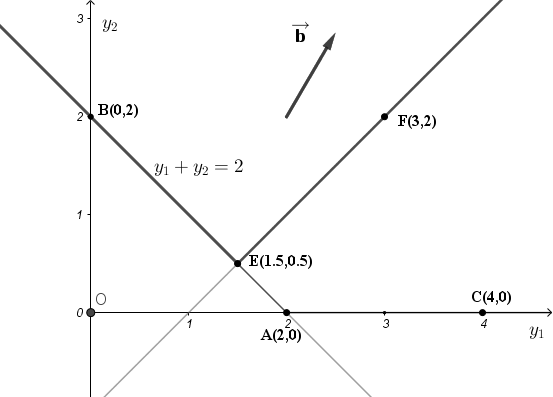}
			\caption{Cone $C_2$}
			\label{fig2}
		\end{figure}

		In the left column, the base variable indicates that plane $(x_2)$ replaces face $(y_2)$, In Figure \ref{fig2}, the vertex of
		cone $C_2$ is $V_2 =(1.5, 0.5)$, and the new cone is the angle $\angle BEF$, where $B=(0, 2)$, $E=(1.5, 0.5)$ and $QG=(3, 2)$, two edges are listed in the above table :
	$e'_1 = (0.5, 0.5)$, $e'_2 = (-0.5, 0.5)$.
	
		Since there's no positive number in the $c^{\land}$, $V_2$ has no cutter and then it is the dual optimum $y^*=(y^*_1, y^*_2)=(1.5, 0.5)$; According to Simplex algorithm the prime optimal point is $x^*=(x^*_1, x^*_2) =(1.5, 0.5)$; The optimal value is 3.5. 
		
		\textbf{The example is end}.
		\section{The $t$-value table and highest cutting method}
		The theory of cone cutting provides a new perspective for Simplex methods, an important gift is the $t$-value table introduced as follows:
		\subsection{The $t$-value table}
		\textbf{Definition 3.1} Given a Simplex table \ref{Tab2.1}, the $t$-value table is the matrix $(t)=(t_{ij})_{m\times (m+n)}$.
		\begin{equation}\label{3.1}
			t_{ij} = c^{\land}_{j} / \tau_{ij}  
		\end{equation}
		
		Where $t_{ij}$ is vacancy if $\tau_{ij}=0$, This table is easily constructed, with each row divides $c^{\land}_j$ row item by item. $Q_{ij}$ is called virtual if it locates on the inverse ray of an edge, for short, say it locates a virtual edge; when $t_{ij}$ is vacancy, $L_i$ parallels or lie in $(\tau_{j})$. According to \ref{2.8}:
		\begin{equation}\label{3.2}
		t_{ij} = e^{*}_{i}A_j  
		\end{equation}
		
		Where $t_{ij}\tau_{ij}=c^{\land}_{j}.$ By means of the $t$-value table, we can quickly determine the location information around intersected points on each edge of a cone. The usage is not only taken for edge, but on more rays starting from any point $P$ along any direction $d$. Generally: Given a ray $R : y=Q+td (t\geq 0)$, writing
		three vectors as follows:
		\begin{equation}\label{3.3}
		\begin{aligned}
		M=(M_1,...,M_{m+n}): M_j=c^+-QA^+_j (j=1,...,n+m);\\
		D=(D_1,...,D_{m+n}): D=dA^+_j (j=1,...,n+m);\\
		t=(t_1,...,t_{m+n}): t_j=M_j/D_j(j=1,...,n+m);
		\end{aligned}
	\end{equation}
	
		We call that $Molecular Denominator$ and $t$-value vectors respectively, and they form a matrix,
		called the Tri-rows.
		\subsection{Highest algorithm}
		We would like to draw a comparison between the classical Simplex algorithm in solving the primary LP problem verse the Cone-Cutting algorithm tackling the dual LP problem. In the primary LP problem solving
		$max\{cx|Ax\leq b\}$, the classical Simplex algorithm starts off with a feasible solution (notice that when $b>0$, the origin $x=0$ is a feasible solution). Then the algorithm identifies a non-basic variable $x_{j^*}$ that can improve the objective function firstly, and then finds a basic variable $x_{i^*}$ that can replace and maintain feasibility. 
		
		For the dual LP problem of $min\{yb|yA\geq c, y\geq 0\}$, the vertex of the coordinate cone (i.e. the origin $O$) is not necessary a feasible point. If it is, then the minimum solution $Vb=0$ is the optimum solution. Otherise, it sacrifices the minimum value but move to another cone that is less infeasible in each iteration. The Cone-Cutting algorithm is to choose a cutting facet $(\tau_{j^*})$ or a cutter that rejects the vertex to cut an existing cone; and then to identify a cut point $\tau_{i^*,j^*} $ which becomes the vertex of a new cone. In moving to the new cone, the facet $(\tau_{i^*})$ is deleted. The new objective value or target value is higher than the previous one, but number of facets that rejects the new vertex is non-increasing. 
		In both cases, it involves identifying $j^*$ first, and then uses \ref{2.12} to compute $i^*$.
		
		The key problem of Simplex method is how to choose the base variable, i.e., how to choose the cutter. The principle of deepest algorithm is to maximize the cutting degree $c^{\land}_j$:
		\begin{equation}\label{3.4}
		j^*= \mathop{\arg\max}_{j}\{c^{\land}_j|c^{\land}_j>0\}.
		\end{equation}
		
		However, this choice may not be the most effective. A measure of effectiveness should be the height difference $(V_{(t)}b-V_{(t-1)}b)$ between the $(t-1)$ iteration and the $t$ iteration. When the norm vector  $\tau_{j^*}$ of cutter deviates from the target vector $b$, the cutting will have little changing on the height of cone vertex. In extreme cases, when $\tau_{j^*}$ is perpendicular to $b$, the height difference is zero. The Klee–Minty anti-example was designed with this in mind. We quote its simplest form as follows:
	\begin{center}
		Table 3.1 Klee-Minty example ($m=3$)
		\begin{equation}\label{Table 3.1}
		\begin{tabular}{|c|ccccccc|c|}
		\hline 
		$	\underline{B}$	& $x_1$ & $x_2$ & $x_3$ & $x_{...}$ & $y_1$ & $y_2$ & $y_3$ & $s$\\ 
		\hline 
		$y_1$	&  1& 0 & 0 & 1 & 0 &9  &  & 1 \\ 
		$y_2$	& 20 &  1& 0 & 0 &1  & 0 &  & 100 \\ 
		
		$y_3$	& 200 & 20 & 1 & 0 & 0 & 1 &  & 10000 \\ 
		\hline 
		$c^{\land}$	& 100 &10  & 1 & 0 & 0 & 0 &  & 0 \\ 
		\hline 
		\end{tabular} \notag
		\end{equation}
	\end{center}
	
		The traditional algorithm, according to (3.4), should choose the first one as the basis, so that the optimal solution can be obtained by $7=2^3 - 1$ pivoting. By extension, the optimal solution can be obtained through $2^m - 1$ pivoting, so it is determined that the deepest algorithm is not got in polynomial time. Now, we write the $t$-value table of Table 3.1 as follows.\\
	\begin{center}
		Table 3.2 The $t$-value table of last table
		\begin{equation}\label{Table3.2}
		\begin{tabular}{|c|ccccccc|c|}
		\hline 
		$	\underline{B}$	& $x_1$ & $x_2$ & $x_3$ &  $y_1$ & $y_2$ & $y_3$ & $s$\\ 
		\hline 
		$y_1$	&  100&  &  & 0  &  &  & 1 \\ 
		$y_2$	& 5 &  10&  &  &0  &  & 100 \\ 
		
		$y_3$	& 0.5 & 0.5 & 1 & && 0 & 10000 \\ 
		\hline 
		$c^{\land}$	& 100 &10  & 1 & 0 & 0 & 0  & 0 \\ 
		\hline 
		\end{tabular} \notag
		\end{equation}
	\end{center}

Where each value $t_{ij}$ reflects the location information of $Q_{ij}$ intersected by plane $(\tau_j)$ and edge $L_i: t_{ij}>0$ implys that $Q_{ij}$ is real, which locates on $\underline{L_i}$; while $t_{ij}<0$ implying that $Q_{ij}$ is virtual, which is in the inversed ray of $\underline{L_i}$. What we care about is the height of these intersections. We have that $Q_{ij}b=Vb+t_{ij}e_ib$, since $Q_{ij}=V+t_{ij}e_i$, and we have that $s_i = e_i b,$ which was emphasized in Theorem \ref{theo2.1}, so that $t_{ij} s_i = Q_{ij}b-Vb$. The right side of the equation is called the relative height of $Q_{ij}$ with respect to the vertex $V$. Times slop $s_i$ to all real $t-$values $t_{ij}$ in $i-$th row, we get the following tableau, called the intersected relative heights table.
	\begin{center}
	Table 3.3 Intersected relative heights table
	\begin{equation}\label{Table3.3}
	\begin{tabular}{|c|ccccccc|c|}
	\hline 
	$	\underline{B}$	& $x_1$ & $x_2$ & $x_3$ &  $y_1$ & $y_2$ & $y_3$ & $s$\\ 
	\hline 
	$y_1$	&  100&  &  & 0  &  &  & 1 \\ 
	$y_2$	& 500 &  100&  &  &0  &  & 100 \\ 
	
	$y_3$	& 5000 & 5000 & 10000 & && 0 & 10000 \\ 
	\hline 
	$c^{\land}$	& 100 &1000 & 0 & 0 & 0 & 0  & 0 \\ 
	\hline 
	\end{tabular} \notag
	\end{equation}
\end{center}

What does the last row of this table represent? For any $j\in J$, let
\begin{equation}\label{3.5}
h_j=min_j\{t_{ij}s_i|t_{ij}>0\}
\end{equation}
which is called the lowest positive height difference.
The highest cut point principle to select a cutter $(\tau_{j^*})$ is:
\begin{equation}\label{3.6}
j^*=\mathop{\arg\max}_{j}\{h_i|c^{\land}_j>0\}
\end{equation}

If for a particular $j\in J$, $t_{i,j}\leq 0$ hold for all $i$, then \ref{3.5} is meaningless. However, if $c^{\land}_j>0$ in this case, then cutter $(\tau_{j})$ will cut the whole acceptable region of the cone out, resulting an empty dual feasible region. This is explained clearly in section 4. Assuming the dual feasible region is non-emplty, \ref{3.6} can always produce a proper $j^*$ solution.

After obtaining $j^*$, \ref{3.5} can be used to find the facet $i^*$ that will leave the basic facet set.
\begin{equation}\label{3.7}
i^*=\mathop{\arg\min}_{i}\{t_{ij}s_i|t_{ij}>0\}
\end{equation}

Because $C$ is a strictly regular cone, all the edge slopes $s_i>0$, so $\tau_{i,j}>0$ when and only when $t_{i,j}>0$; since $c^{\land}_{j^*}$ is not dependent on $i$, we have
\begin{equation}
i^*=\mathop{\arg\min}_{i}\{t_{i,j*}s_i|t_{i,j*}>0\}=\mathop{\arg\min}_{i}\{c^{\land}_js_i|t_{i,j*}>0\}=\mathop{\arg\min}_{i}\{s_i/\tau_{i,j*}|t_{i,j*}>0\}=i(j^*),\notag
\end{equation}

This is identical to \ref{2.12} in computing the leaving basis in Simplex. Hence in the highest algorithm, when the cuter $j^*$ is selected, it uses the same method as Simplex to determine the pivot point $(i^*, j^*)$. In the Klee-Minty counter-example, after determining the pivot point $(i^*, j^*)=(3,3)$, it performs the pivot operation on \ref{Table 3.1} to obtain the following tableau:
			\begin{center}
						Table 3.4 Pivoting from Table 3.1
					\begin{equation}\label{Tab3.4}
					\begin{tabular}{|c|ccccccc|c|}
					\hline 
					$	\underline{B}$	& $x_1$ & $x_2$ & $x_3$ &  $y_1$ & $y_2$ & $y_3$ & $s$\\ 
					\hline 
					$y_1$	&  1& 0 & 0 & 1  & 0 & 9 & 1 \\ 
					$y_2$	& 20 &  1& 0 & 0 &1 & 0 & 100 \\ 
					
					$y_3$	& 200 & 20 & 1 & 0&0& 1 & 10000 \\ 
					\hline 
					$c^{\land}$	& -100 &-10  & 0 & 0 & 0 & -1  & -10000 \\ 
					\hline 
					\end{tabular} \notag
					\end{equation}
				\end{center}
	
				Since that there is no positive number occurring in $c^{\land}$, we get the dual optimal point is $y^*=(0, 0, 1)$, with optimal height $h^{*}=10000$. The highest algorithm gets the solution of Klee-Mimty counter-example by only one time of pivoting.
				
			In our papers \cite{13}, there was a gravity sliding algorithm, which is not a method of cone-cutting, but a method of feasible point falling and sliding along the wall of feasible region. The corresponding algorithm in cone-cutting can be described as follows:
			
			The steepest principle for cutter selection:
			\begin{equation}\label{key}
			j^*=\mathop{\arg\max}_{j}\{s_{i(j)}\}.
			\end{equation} 

The steepest principle can be used to solve the Klee-Minty counter example by one time of cone-cutting also, but we could not say that steepest is better than deepest principle absolutely. We can get counter example showing that deepest can be better than steepest principle. However, The highest principle is always better than deepest and steepest principles both in a large probability. 

%
%
%
\section{Column elimination}
				\subsection{Basic theorem of column elimination}
				In cone cutting, an expelled face can cut back into base cone frequently, this is the main obstacle for the searching of strong polynomial solution. The column elimination theorem proposed by Y. Ye in article \cite{9} has important significance, and we have advanced work stated in this section.\\
				\textbf{Definition 4.1} A constraint plane ($\tau_{j}$) is called an redundant plane if it does not intersect the feasible region $D$ and it accepts $D$; ($\tau_{j}$) is called a golden plane if it includes an optimal point; ($\tau_{j}$) is called a defective plane if it neither the redundant, nor the golden plane. A cone is called a golden cone if all its faces are golden.\\ 
				\textbf{Definition 4.2} A constraint plane ($\tau_{j}$) is called eliminatable if it is not golden.
				
				 How can the eliminatable planes be eliminated? \\
				\textbf{Theorem 4.1} (New elimination theorem) If $\tau_{j}\geq0$ and $c^{\land}_j<0$, then ($\tau_{j}$) is eliminatable; If $\tau_{j}\leq0$ and $c^{\land}_j>0$, then ($\tau_{j}$)  cuts off all feasible points. \\
				\textbf{Proof} According to \ref{3.1}, If $\tau_{j}\geq0$ (i.e., for i=1,...,m, $\tau_{ij}\geq 0$) and $c^{\land}_j\leq0$, then for i=1,...,m, we have that $\tau_{ij}\leq0$. For those index $i$ who makes that $\tau_{ij}=0$, $\tau_{ij}$ parallels to and accepts $L_i$ since $c^{\land}_j\leq0$;
				The rest situation is that all intersected parameters are negative: $\tau_{ij}\leq0$ . It means that all intersected
				points are virtual. Therefore, the whole cone $C$ is located at one side of the plane ($\tau_{j}$). Since $c^{\land}_j<0$, ($\tau_{j}$)	accepts the whole cone $C$ but does not contain any optimal point. So that ($\tau_{j}$) is eliminatable.
				
				This theorem is similar to the Ye-column elimination theorem, but derived from different basis, and more intuitive.
				
				The condition of elimination theorem is that all coefficients in the norm vector must have a same symbol: non-negative or non-positive, we call such a property as symbol consistency. How to take a plane having symbol consistency? We need to introduce the following concepts.
				\subsection{Implantation of the reverse horizontal} 
				\textbf{Definition 4.3} A normal cone called a strictly normal cone, if for $i=1,…,m, s_i>0$. A non-strict normal cone is a normal cone with horizontal edges.
				
				The strictly normal cone is unbounded on top, in order to get a bounded enclosure, we put a lid on it.\\
				\textbf{Definition 4.4} For a given positive number $u$, the $u$-reverse horizontal lid is a plane with equation in the original table $(T_o)$:
				\begin{equation}\label{4.1}
					-b_1y_1-...-b_my_m=-u
				\end{equation}
			
				Or, in the table \ref{Tab2.2}:
			\begin{equation}\label{4.1'}
			-s_1y_1-...-s_my_m=-u \tag{4.1'}
			\end{equation}
			
		Which is called the reverse horizontal lit of a strict normal cone $C$, Note that the norm vector of the lit is $-b$, or $-s=-B^{-1}b$.\\
	\textbf{Definition 4.5} Suppose that the $C$ is a strictly normal cone, denote that
			\begin{equation}\label{4.2}
				C_u=\{y\in C|h(y)\leq u\},D_u=\{y\in D|h(y)=u\}.
			\end{equation}
			
			$C_u$ is called the $u$-frustum of $C$ and $D_u$ is the $u$-feasible section of $C$.\\
			\textbf{Proposition 4.1} Given a Simplex table \ref{Tab2.2}, as long as there is a feasible point whose height is equal to or lower than $u$, then, without affecting the solution of the programming problem, $u$-reverse horizontal plane can be added to the table and the feasible region $D$ can be substituted by $D_u$.\\
			\textbf{Proof} The set $D_u$ is formed by deleting points higher than $u$ from $D$. Let $y^*$ be a lowest point, because it is the lowest point in $D$, and its height will not be higher than the height of any feasible point, so $y^*$ will not be higher than $u$, and then will not be deleted by the $u$-reverse horizontal plane. Therefore, adding this lit to the table will not change the acquisition of the lowest point, and $D_u$ can be used to replace $D$.  \; \textbf{End}
			
			The proposition benefits us that a defective plane for $D$ can be changed to a redundant plane for $D_u$, the later one is easer to be eliminated than the former.\\
			\textbf{Definition 4.6} Given a Simplex table (\ref{Tab2.2}), add the reverse horizontal plane column with a column index, after the $n$-th column, the expanded table is called an expanded table, denoted as ($T^+$).
			\subsection{Symbol consistency theorem}
			\textbf{Theorem 4.2}\label{th4.2} (Symbol consistency coefficient theorem) Given a strictly normal cone plus an $u$-lit written in the expanded table ($T^+$). Doing pivoting at $(i^*,\Delta)$ according to following $t^*$-values.
			\begin{equation}\label{4.3}
				i^*=\mathop{\arg\max}_{i}\{s_i/\tau_{ij}|\tau_{ij}<0\} (when \; c^{\land}_j<0);
			\end{equation}
			\begin{equation}\label{4.4}
				i^*=\mathop{\arg\min}_{i}\{s_i/\tau_{ij}|\tau_{ij}>0\} (when\; c^{\land}_j>0).
			\end{equation}
		Then, the norm vector of $\tau_{j}$ can become non-negative by \ref{4.3} and non-positive by \ref{4.4}.\\
		\textbf{Proof} Note that the norm coefficients of $(x^\Delta)$ is that $\tau_{i\Delta}=-s_i (i=1,...,m)$, After pivoting, we have that
		\begin{center}
			$\tau'_{ij}=\tau_{ij}-\tau_{i*j}\tau_{i\Delta}/\tau_{i^*\Delta}=\tau_{ij}-\tau_{i*j}(-s_i)/(-s_{i^*})=\tau_{ij}-\tau_{i*j}s_i/s_{i^*}$
		\end{center}
		Hence, according to \ref{4.3}, $s_{i^*}/\tau_{i*j}\geq s_i/\tau_{ij}$, i.e., $\tau_{ij}/s_i\geq\tau_{i*j}/s_{i^*}$, while $s_i>0$, then we have that
		\begin{center}
		$\tau'_{ij}=s_i(\tau_{ij}/s_i-\tau_{i*j}/s_{i^*})\geq0$
		\end{center}
		The new coefficents become non-negative by \ref{4.3}. Similar proof can be got for the rest.  \; \textbf{End}\\
		\textbf{Example 4.1} Given the extended table $(T^+)$, try to sign the coefficients of plane $(x_5)$ .
		\begin{center}
				Table 4.1 extended table $(T^+)$
			\begin{equation}\label{Table4.1}
			\begin{tabular}{|c|cccccccccc|c|}
			\hline 
			$	\underline{B}$	& $x_1$ & $x_2$ & $x_3$ & $x_4$&$x_5$&$x^{\Delta}$& $y_1$ & $y_2$ & $y_3$ &$y_4$& $s$\\ 
			\hline 
			$y_1$	&  -2& -2 & 1 & 1  & 0 & -1 & 1 &0&0&0&1\\ 
			$y_2$	& -1 &  4& -1 & 2 &-1 &-2& 0 & 1&0&0&2 \\ 	
			$y_3$	& 3 & 2 & 0 & -1&1& -6 & 0 &0&1&0&6 \\ 
			$y_4$	& 1 & -1 & 3 & 0&1& -3 &0 &0&0&1&3 \\ 
				\hline 
			$c^{\land}$	& -1&2  & 1 & 2 & -2 & -u  & 0&0&0&0&0 \\ 
			\hline 
			\end{tabular} \notag
			\end{equation}
		\end{center}		
		\textbf{Solve} Since that $c^{\land}_5=-2<0$, according to \ref{4.3} $i^*$=$\mathop{\arg\max}_{i}\{\tau_{ij}/\tau_{ij*}|\tau_{ij}<0\} =2$,
		Doing pivoting at $(2,\Delta)=(2, 6)$, get new table as follows:
	\begin{center}
		Table 4.2 ($x^\Delta$) enters into the base
		\begin{equation}\label{Table4.2}
		\begin{tabular}{|c|cccccccccc|c|}
		\hline 
		$	\underline{B}$	& $x_1$ & $x_2$ & $x_3$ & $x_4$&$x_5$&$x^{\Delta}$& $y_1$ & $y_2$ & $y_3$ &$y_4$& $s$\\ 
		\hline 
		$y_1$	&  -3/2&-4& 3/2 & 3/2 & 1/2  & 0 & 1 & -1/2 &0&0&0\\ 
		$x^{\Delta}$	& 1/2 &  -2& 1/2 & 1/2 &1/2 &1& 0 & -1/2&0&0&-1 \\ 	
		$y_3$	& 6 & -10 & 3 & 2&4& 0 & 0 &-3&1&0&0 \\ 
		$y_4$	& 5/2 & -7 & 9/2 & 3/2&5/2& 0 &0 &3/2&0&1&0 \\ 
		\hline 
		$c^{\land}$	& \multicolumn{9}{c}{-2+u/2 \; -u}& &-u \\ 
		\hline 
		\end{tabular} \notag
		\end{equation}
	\end{center}	

		This pivoting is taken at a negative grid, which is a non-classical Simplex. The approach is still to
		normalize the pivoting grid, times$-1$ to second row, and then continue to operate.
		
		Note that a non-classic pivoting gets no longer a normal cone, the second edge slope of the edges is negative, The transformation brings the norm coefficients of plane $(x_5)$ in consistency: $\tau_{5}=(1/2, 1/2, 4, 5/2)^T$.
		
		The intuitive idea of Theorem 4.2 is: Treating $C_u$ as a cone with a flat top cover at height $u$ , there are $m$ vertexes $V_1(u),...,V_m(u)$ on the cover. We refer these as edge-vertices to distinguish it from the vertex of a cone. $C_u$ becomes an enclosed polyhedron. Let $u$ began to rise from $h(V)$, The whole body of the intersected
		cone $C_u$ was located in the lower side of any given plane $(\tau_{j})$, Once a vertex $V_{i^*(u)}$ contacts with the $\tau_{j}$,
		firstly at $u=u^*$, the plane then gets symbols consistency there: Let $C_{i*}$ be the cone, which is got by doing pivot at $(i^*,\Delta)$, i.e., which is the original cone $C$ cut by $(x^{\Delta})$. It is obvious that $V_{i^*(u)}$ is the vertex of $C_{i^*}$. If $u=u^* - \epsilon$ then $(\tau_{j})$ cuts $C_{i^*}$ on real edges consistently; If $u=u^* + \epsilon$, then $(\tau_{j})$ cuts $C_{i^*}$ on virtual edges
		consistently. The vertex $V_{i^*(u)}$ is so important for us, How to get it? It is the lowest real intersection of ($\tau_{j}$) with the older cone $C$. This is what mean in \ref{4.3} and \ref{4.4}.\\
		\textbf{Definition 4.8} The $i$-th edge is called the fishing edge with respect to $(x_j)$; the intersected point $V_{i^*}$ is
		called its critical point; $u^*=h(V_{i^*})$ is called the critical height of ($x_j$). $u$ is called a feasible height if there is a feasible point $Q$ with $h(Q)\leq u$.
		
	The above discussion leads to the following theorem:\\
		\textbf{Proposition 4.2} Every constraint plane ($x_j$) has at least one fishing edge; each fishing edge has one and only one critical point $V_{i^*}$:
		\begin{equation}\label{4.5}
		V_{i^*}=V+c^{\land}_je_{i^*}/\tau_{i*j}
		\end{equation}
		
	with creticle height $u_{j}^*$
\begin{equation}\label{4.6}
u_{j}^*=h(V)+c^{\land}_js_{i^*}/\tau_{i*j}
\end{equation}

The critical point $V_{i^*}$ of plane ($x_j$) implies that ($x_j$) could be deletable provided its critical height $u_{j}^*$ is a feasible height;
otherwise, if all feasible points exist above critical height only, then no meaning for critical point again. So far, we have following theorem:\\
\textbf{Theorem 4.3} A plane ($x_j$) with $c^{\land}_j<0$ can be deleted if and only if there is a dual feasible point being not higher than its criticle point; A plane ($x_j$) with $c^{\land}_j>0$ can cut off the whole dual feasible region $D$ , if and only if there is a dual feasible point being not higher than its critical point.

The theorem tells us that we can't guarantee that all non-golden columns can be deleted immediately; column elimination relies on the falling down of feasible
height, we need to wait about feasible falling stated in Section 6.\\
	\textbf{Example 4.2}  Critical height records.
	\begin{center}
		Table 4.4 Records of critical heights
		\begin{equation}\label{Table4.2}
		\begin{tabular}{|c|ccccccccccc|c|}
		\hline 
		$	\underline{B}$	& $x_1$ & $x_2$ & $x_3$ & $x_4$&$x_5$&$x_6$& $y_1$ & $y_2$ & $y_3$ &$y_4$& $y_5$&$s$\\ 
		\hline 
		$x_1$	&  1&3.5& 0.5 & 0 & 0  & 0 & 0.5 & -0.5 &0&0&0&1.5\\ 
		$x_4$	& 0 &  2.5& 0.5 & 1 &0 &1& -0.5 & 1.5&0&1&0&1.5 \\ 	
		$y_3$	& 0 & -2 & -1 & 0&0& -1 & 1 &-2&1&-2&0&2 \\ 
		$x_5$	& 0 & 1.5 & 0.5 & 0&1& 0 &-0.5 &0.5&0&1&0&0.5 \\ 
		$y_5$   & 0 & 0.5 & -1.5 & 0&0& 1 &-0.5 &0.5&0&0&1&0.5 \\ 
		$c^{\land}$   & 0 & -19 & -2 & 0&0& 0 &-3 &1& 0& 0&0 &3.5 \\ 
			\hline 
		 &\multicolumn{2}{r}{-30}&\multicolumn{4}{l}{-15}&\multicolumn{4}{l}{-14}&&-u \\ 
		\hline 
		\end{tabular} \notag
		\end{equation}
	\end{center}	

	In theory, we need to use expanded table, but we can use table \ref{Tab2.2} directly according to proposition 4.2.
	
	The cutting degrees of ($x_2$), ($x_3$) and ($y_1$) are negative, and their critical heights can be calculated
	according to \ref{4.6} respectively:
	
	Fishing index of $(x_2)$: $i^*=3, u_2*=(V)+c^{\land}_2s_3/\tau_{32}=30;$ 
		
	Fishing index of $(x_3)$: $i^*=3, u_3*=(V)+c^{\land}_3s_3/\tau_{33}=15;$
	
	Fishing index of $(y_1)$: $i^*=3, u_6*=(V)+c^{\land}_6s_4/\tau_{46}=14;$
	
	\textbf{End of Example. 4.1}
	
	We only introduce how to record eliminatable information for negative cutting degree plane based on \ref{4.3}, but do not introduce how to treat cut off situation for positive cutting degree planes based on \ref{4.4}.
	Why? Because, an Assumption for strong polynomial problem in the paper is that: The dual feasible region
	$D$ is not empty. So there is no necessary to consider the cut off situation. It is not no such desire, it is no
	enough ability.
	
	The same plane can have positive, negative or zero cutting degree in different tables. Therefore, the
	critical height records of a plane may appear and disappear, and may contradicted each other, but the plane
	can be eliminated as long as one record is feasible.
	\section{Ray feasible interval and non-strict normal cone treatment}
	The horizontal upper cover of the cone requires that the cone is strictly normal, no edges with a zero slope. Otherwise, the horizontal intersected cone is not a bounded closed convex set, and the pivoting grid will stays on the horizontal edge, which causes the taboo of equal height cycle in the traditional Simplex algorithm. In fact, the appearance of the horizontal edge is exactly the information what we look for firstly, if there is a feasible point on the horizontal edge, then all the points in the feasible interval are the dual optimal points.
	\subsection{Looking for feasible interval on a ray}
	
	Given the Simplex table \ref{Tab2.1}, rays are extracted from point $P$ along direction $d$, and we are looking for the feasible internal on the ray.\\
	\textbf{Theorem 5.1} Judgement of feasible interval on a ray Given a Simplex table \ref{Tab2.1}. According to Definition 3.1:
	
	\textbf{F1)} If there is a $D_j=0$ with $M_j>0$, Then there is no feasible point on the ray;
	\textbf{F2)} If $M_j\leq0$ whenever $D_j=0$, then set that
	\begin{equation}\label{key}
		t_a=\mathop{\max}_j\{t_j|D_j>0\};
	\end{equation}
\begin{equation}\label{key}
t_b=\mathop{\min}_j\{t_j|D_j<0\};
\end{equation}
	
	If $t_a\leq t_b$ then there is a feasible interval $[Q_a, Q_b]$ on the ray:
	\begin{equation}\label{key}
	Q_a=P+t_{a}d,\; Q_b=P+t_{b}d,
	\end{equation}
	
	Proof is obvious.
	
	It is wrong that accepting a feasible interval by means of condition \textbf{F2)} only. We must check if the
	condition \textbf{F1)} is satisfied in the first. When $P=V$ is the vertex of cone C with respect to Table $(T)$ and $d=e_i^*$ is the direction of edge $\underline{L_i}$ of $C$
	then the Molecular vector 
	$M$ is the lowest row $c^{\land}$ in table $(T)$, and the Denominator vector $D$ is the $i$-th
	row in table $(T)$.
	
	It is wrong that calculating $D_j=d\tau_{j}$ in stead of $D_j=dA_{j}^+$. Since that $\tau_{j}$ is the variant of $A^+_j$, it is very
	easy to mix the usages of them. When you calculate the denominator $D$, please turn back to the original
	table \ref{Tab2.1}.
	\subsection{Treatment of non-strict normal cones}
	\textbf{Example 5.1} Discover a horizontal edge from Tables 5.1.
	\begin{center}
		Table 5.1 Original table $T_o$
		\begin{equation}\label{Table5.1}
		\begin{tabular}{|c|ccccccccccc|c|}
		\hline 
		$	\underline{B}$	& $x_1$ & $x_2$ & $x_3$ & $x_4$&$x_5$&$x_6$& $y_1$ & $y_2$ & $y_3$ &$y_4$& $y_5$&$\Delta h$\\ 
		\hline 
		$y_1$	&  2&-2& 1 & 1 & 0  & 1 & 1 & 0 &0&0&0&4\\ 
		$y_2$	& 0 &  1& 0 & 1 &-1 &1& 0 & 1&0&0&0&1 \\ 	
		$y_3$	& 0 & 2 & 0 & 1&1& 0 & 0 &0&1&0&0&4 \\ 
		$y_4$	& 1 & 1 & 1 & 0&1& 0 &0 &0&0&1&0&2 \\ 
		$y_5$   & 1 & 0 & -1 & 0&0& 1 &0 &0&0&0&1&6 \\ 
		\hline
		$\sigma$   & 6 & -20 & 1 & -2&2& 0 &0 &0& 0& 0&0 &0 \\ 
		\hline 
		\end{tabular} \notag
		\end{equation}
	\end{center}	
	
	This is tableau represents a strictly regular cone, but it has potential horizontal edges. Generally speaking, if the exit basic variable $i^*=i(j)$ of the $j^{th}$ column is not unique, a horizontal edge must appear when the pivot operation is performed on $(i^*,j)$. In the first column of this table, $i^*=i(j)=1$ or $4$. If (4,1) is used as the pivot point, the first edge of the new tableau after the transformation is a horizontal edge. This can be observed as the edge slope $s_1$ is zero. Hence $\underline{L1}$ is a horizontal edge.
	
	Once the horizontal edge appears, it is necessary to find the feasible segment on it to determine the optimal solution. According to Section 5.1, since the ray we are considering is an edge of a cone, the vector $t$ is a row vector of the $t$-value table mentioned in the previous section. Thus, we can write down the three-row ray cut matrix based on the values of the tableau $(T)$ as follows: The numerator vector is the last row of $(T)$, the denominator vector is the first row of $(T)$, corresponding to the horizontal edge; and the $t$ value vector can be readily obtained by doing element-wise division. This is shown below:
	
	It is a criterion
	that if $s_i/\tau_{ij}=s_k/\tau_{kj}$, then taking pivot at $(i,j)$ or at $(k,j)$, there must occur the horizontal
	edge in the next table. Now, do pivot at $(4,1)$, get table 5.2.
	\begin{center}
		Table 5.2 Discover an horizontal edge
		\begin{equation}\label{Table5.1}
		\begin{tabular}{|c|ccccccccccc|c|}
		\hline 
		$	\underline{B}$	& $x_1$ & $x_2$ & $x_3$ & $x_4$&$x_5$&$x_6$& $y_1$ & $y_2$ & $y_3$ &$y_4$& $y_5$&$s$\\ 
		\hline 
		$y_1$	&  0&-4& -1 & 1 & -2  & 1 & 1 & 0 &0&-2&0&0\\ 
		$y_2$	& 0 &  1& 0 & 1 &-1 &1& 0 & 1&0&0&0&1 \\ 	
		$y_3$	& 0 & 2 & 0 & 1&1& 0 & 0 &0&1&0&0&4 \\ 
		$y_4$	& 1 & 1 & 1 & 0&1& 0 &0 &0&0&1&0&2 \\ 
		$y_5$   & 0 & -1 & 2 & 0&-1& 1 &0 &0&0&-1&1&4 \\ 
		\hline
		$\sigma$   & 0 & -26 & -5 & 2&-8& 2 &0 &0& 0& -6&0 &-12 \\ 
		\hline 
		\end{tabular} \notag
		\end{equation}
	\end{center}	
	
	Since that $s_1=0$, the $1$-th edge is a horizontal edge, then looking for the feasible interval on the edge $\underline{L_1}$, according to Theorem 5.1:
	\begin{center}
		\begin{tabular}{cccccccccccc}
			&$x_1$ &$x_2$ &$x_3$ &$x_4$ &$x_5$ &$x_6$ &$y_1$ &$y_2$ &$y_3$ &$y_4$ &$y_5$ \\ 
		M: &0 &-26& -5& 2& -8& 2& 0& 0&  0 &-6& 0\\ 
			D:&  0 & -4&  -1&  1&  -2&  1 & 1 & 0&  0&  -2&  0\\ 
			t:& \multicolumn{2}{r}{6.5}& 5&  2 &4& 2 &0 & & &3 &\\ 
		\end{tabular} 
	\end{center}
	
	F1. No positive cutting degree when $D$-value is zero, satisfied;
	
	F2. $t_a=max_{j}\{t_j|\tau_{1j}>0\}=2<3= min_j\{t_j|\tau_{1j}<0\}=t_b,$
	
	There is a feasible interval on the edge:
	
	$Q_a=(0,0,0,6,0)+2(1,0,0, -2,0)=(2,0,0,2,0)$
	
	$Q_b=(0,0,0,6,0)+3(1,0,0,-2,0)=(3,0,0,0,0)$
	
	Conclusion: The dual optimal points are the interval $[Q_a,Q_b]$, the optimal value is $h(Q_a)=12$.	
	
	Most optimal solutions in big data linear programming are company with horizontal edge problem. According to
	this paper’s idea, the most difficult problem becomes easier task.
	\section{Horizontal feasible central falling}
	The elimination algorithm in Section 4 leaves a task to this section: Identify whether a height is feasible. There needs given a feasible point $F$ in the dual feasible region $D$. The height $u_o=h(F)$ is feasible, then, let the feasible point fall down, the feasible height will be extended down also. In paper \cite{13,14}, the authors have put forward a feasible point falling along the gravity direction $g=-b$ until down to the wall of $D$. then, sliding along the projection of $g$ in blocking planes until down to an optimal point, named the \emph{Gravity Sliding algorithm}. Which realizes the maximal gradient principle on the convex polyhedrons. Despite its advantages, the idea to be updated in this section: When the feasible point hits the wall, instead of sliding in boundary planes, we move the point horizontally to the inner of $D$ such that it falls down along the gravity g continuously, so to maintain as inner optimization.
	\subsection{Horizontal feasible outline}
	Given an expanded table $(T^+)$. Let $Q$ be a given feasible point.
	
	Calculate $u^*=h(Q)=Qb$, assign that $u=u^*-h(V)$ in $(T^+)$, which is relative height of $Q$ with respect
	to the cone vertex $V$.
	
	For $i=1,...,m$,\\
	1) Doing pivoting at $(i,\Delta)$, the $i$-th vertex $V_i$ of the upper cover of $C_{u^*}$ is got.\\
	2) Set $P_i=V_i,\; d_i=Q - P_i$, and find out the feasible interval on the ray $y=P_i+td_i$, Since that the ray
	starts from $P_i$ and directs to $Q$, we are looking for the start point of the feasible interval, which is
	nearer the vertex $V_i$, denote the first point by $F_{\Delta}$.\\
\textbf{Definition 6.1} Denote:
	\begin{equation}\label{key}
		[F_{\Delta}]=\{y|\exists (\lambda_1,...\lambda_m)\geq 0,\lambda_{1}+,...,\lambda_{m}=1;\lambda_1F_1+,...,\lambda_{m}F_m=y\}
	\end{equation}
	
	Which is called the horizontal feasible outline on $(x^{\Delta})$. 
	
		$[F_{\Delta}]$ may not include all feasible points on $(x^{\Delta})$, which is just an outline.\\
	\textbf{Example 6.1} Given a Simplex table \ref{Tab2.2} (the same as Table 5.2), given a feasible point $F=(3.5, 0.2, 0, 0.2, 0)$, try to find out the supportor of the upper cover containing $F$. \\
	Step (0): Calculate the height of $F: h(F)=Fb=F(4,1,4,2,6)^T=14.6$.
	
	Calculate the height of $V=(0,0,0,6,0)$, the vertex of $C$ in Table 5.2: $h(V)=12$;\\
	Assign $u=14.612=2.6$ in the expanded table $(T^+)$:
		\begin{center}
		Table 6.1 Assign $u=2.6$ in $(T^+)$
		\begin{equation}\label{Table6.1}
		\begin{tabular}{|c|cccccccccccc|c|}
		\hline 
		$	\underline{B}$	& $x_1$ & $x_2$ & $x_3$ & $x_4$&$x_5$&$x_6$&$x^{\Delta}$& $y_1$ & $y_2$ & $y_3$ &$y_4$& $y_5$&$s$\\ 
		\hline 
		$y_1$	&  0&-4& -1 & 1 & -2  & 1 &0& 1 & 0 &0&-2&0&0\\ 
		$y_2$	& 0 &  1& 0 & 1 &-1 &1&-1& 0 & 1&0&0&0&1 \\ 	
		$y_3$	& 0 & 2 & 0 & 1&1& 0 & -4&0 &0&1&0&0&4 \\ 
		$y_4$	& 1 & 1 & 1 & 0&1& 0&-2 &0 &0&0&1&0&2 \\ 
		$y_5$   & 0 & -1 & 2 & 0&-1& 1 &-4&0 &0&0&-1&1&4 \\ 
		\hline
		$\sigma$   & 0 & -26 & -5 & 2&-8& 2 &-2.6 &0& 0& 0&-6&0 &-12 \\ 
		\hline 
		\end{tabular} \notag
		\end{equation}
	\end{center}
	
	Note that the norm vector of $(x^\Delta)$ was $\tau_{\Delta}=-b$in the original table, but now, $\tau_{\Delta}=-B^{-1}b=-s$. For $i=1,...,5$, calculating feasible start points on $i$-th ray(from $V_i$ directs $F$), here $i=3$.\\
	Step (1): Doing pivoting at $(3,\Delta)$, get Table 6.2: 
		\begin{center}
			Table 6.2 Doing pivot at $(3,\Delta)$ on table 6.1
		\begin{equation}\label{Table6.2}
		\begin{tabular}{|c|cccccccccccc|c|}
		\hline 
		$	\underline{B}$	& $x_1$ & $x_2$ & $x_3$ & $x_4$&$x_5$&$x_6$&$x^{\Delta}$& $y_1$ & $y_2$ & $y_3$ &$y_4$& $y_5$&$s$\\ 
		\hline 
		$y_1$	&  0 &-4   & -1 & 1 & -2  & 1  &0  & 1 & 0  & 0 &-2 &0 &0\\ 
		$y_2$	&  0 &  0.5 & 0 & 0.75 &-1.25  &1&0 & 0 & 1&-0.25 &0 &0 &0 \\ 	
		$y_3$	&  0 & -0.5  & 0 & -0.25 &-0.25   & 0  & 1 &0 &0 &-0.25 &0 &0  &-1 \\ 
		$x_1$	&  1 & 0     & 1 & -0.5 &0.5   & 0 &0  &0 &0 &-0.5 &1 &0  &0 \\ 
		$y_5$   &  0 & -3    & 2 & -1&-2  & 1 &0 &0  &0 &-1 &-1 &1 &0 \\ 
		\hline
		$c^{\land}$   & 0 & -27.3 & -5 & 1.35&-8.7& 2 &0 &0& 0& 0&-0.65 &-6 &-14.6 \\ 
		\hline 
		\end{tabular} \notag
		\end{equation}
	\end{center}
	
	Since the number on $(3,\Delta)$ is negative, this is a non-standard pivoting, by multiplying the third row by $-4$ and so on, operations will continue as usual.
	
	The vertex of the new cone is $V_3=(0,0,0.65,6,0)$.\\
	Step (2): Set $d_3=F-V_3=(3.5,0.2, -0.65, -5.8, 0)$. Searching the feasible interval on the ray $y=V_3+td_3$ on feasible interval, get the start point of feasible interval: $F_3=(1.9, 0.1, 0.3, 2.85, 0)$.
	
	Since the number in grid $(1,\Delta)$ is zero, we can’t do pivoting. we can take the direction $d=e_1$ and calculating the first point of feasible interval on the ray and calculating the first point of feasible interval on the ray $y=V_1+t_{e_1}$ and get the result. Finally we can get the horizontal feasible
	outline expanded by the following 5 vertices:
	\begin{center}
		$F_1 =(3.6, 0.2, 0, 0)$;\\
	$F_2=(0, 2.6 , 0, 6 , 0)$;\\
	$F_4=(1.9 , 0.1 , 0 , 3.45,0)$;\\
	$F_5=(1.8, 0.2, 0, 2.7, 0.3)$;\\
	$\underline{F}=(1.84, 0.64 , 0.06, 3, 0.06)$;
	\end{center}
	\subsection{Feasible point central falling}
	\textbf{Definition 6.2} Let $F$ be a feasible point, we call $y=F+tg=F-tb(t\geq0)$ the falling ray from $F$. The end point of the feasible interval on the ray is called the foot of $F$, denoted as $F\downarrow$, inversely, we call $F$ the lift of $F\downarrow$.
	
	We are going to set $F$ as the center of horizontal feasible outline the feasible central falling algorithm can be given as follows:
	
	\textbf{Input} A strictly normal simplex table $T_o$ with $b\geq0$ and a feasible point $F_o$;

	Step 1: $F:= F_o $; calculating the certer $\underline{F}$ of $[F_\Delta]$; 
	
	Step 2: Calculating the landing point $\underline{F}\downarrow$ from the falling ray $y=\underline{F}-tb$ ($t\geq0$).
	
	Step 3: Selecting a cutter passing through $\underline{F}\downarrow$ and doing cone cutting; If get the optimal point, then stop; Else, go back step 1.

	\textbf{In Example 5.1} we have that 
	\begin{center}
	$\underline{F}=(F_1+F_2+F_3+F_4+F_5)/5=(1.84, 0.64, 0.07, 3, 0.06)$.\\
	\end{center}

And set $d=-b$, searching the feasible interval on the falling ray $y= \underline{F}+td=\underline{F}-tb$:
	\begin{center}
		\begin{tabular}{cccccccccccc}
			&$x_1$ &$x_2$ &$x_3$ &$x_4$ &$x_5$ &$x_6$ &$y_1$ &$y_2$ &$y_3$ &$y_4$ &$y_5$ \\ 
			M: &-0.9 &-20& -3.8& -0.6& -4.5& -0.6& -1.9& -0.6&  -0.1 &-3& -0.1\\ 
			D:&  -16 & -3&  0&  -9&  -5&  -11 & -4 & -1& -40&  -2&  -6\\ 
			t:& 0.06&6.67& & 0.07&  0.9 &0.05& 0.48 &0.6 &0.03 &1.5 &0.01 \\ 
		\end{tabular} 
	\end{center}

Since the starting point is feasible, the $t$-value of landing point is the minimum of positive numbers in the $t$-
vactor, which is 0.01 now, and the landing point is
\begin{equation}\
\underline{F}\downarrow = \underline{F}-0.01b=(1.8, 0.63, 0.03, 2.98,0)  \notag
\end{equation}
We have that
\begin{equation}\
h(\underline{F}\downarrow)= 13.9, h(\underline{F})-h(\underline{F}\downarrow)=14.6-13.9=0.7.  \notag
\end{equation}

While $h(\underline{F})-h(V)=14.6-12=2.6.$

The effective rate of the central falling is
$0.7/2.6 =27\%$, it is not so good in the example, we can return to the Step
1 and get more lower feasible height, but we do not state in detail here.

Review the Example 4.2, there were 3 critical points: The criticle point of $(x_2)$ with height 30 which is not lower
than the feasible height 13.9, so that, the plane $(x_2)$ can be deleted Similarly, since $(x_3)$ has a critical point with height
15 and $(y_1)$ has a critical point with height 14, they can be deleted both.

Be careful, we need not delete those columns from the simplex tableau, we just delete their indices from the index-set $J$. Whenever the number $|J|=m$, it means that the golden cone has been got.

Significance of feasible Central Falling method: The falling velocity has nothing to rely on the dimension. For a
two dimensional feasible region, falling from a feasible center requires not hitting the lowest point but only falling close
enough to the lowest point, which is a simple problem without complexity. We can start cone cutting from landing
position as close as we want to the optimal point, that's what we are fascinated.
\subsection{Variable weight combinations}
	\textbf{Definition 6.3} Let $w=(w_1,...,w_m)$ be a vector with coefficients non-negative and summed 1, denote that
	\begin{equation}\label{6.1}
		F^*=w_1F_1+...+w_mF_m
	\end{equation}
	
	Called the convex combination of $F_1,...,F_m$, $w$ is called the weight vector. When $w$ is variable, $Q^*$ is called a variable weight combination.\\
	\textbf{Theorem 6.1} Homomorphic Theorem Let $t_k$ be the falling $t$-vector from $F_k$, $t^*$ is falling $t$ -vector from
	$F^*=w_1F_1+...w_mF_m$, we have that
	\begin{equation}\label{6.2}
	t^*=w_1t_1+...+w_mt_m
	\end{equation}
	\textbf{Proof} The denominator $D=F^*A^+=(w_1F_1+...+w_mF_m)A^+=w_1F_1A^++...+w_mF_mA^+$\\
	Since $w_1+...+w_m=1$, we have that
	\begin{center}
	$M=(w_1+...+w_m) c^+=w_1c^++...+w_mc^+$,
	$M=c^+-Q^*A^+=w_1(c^+-Q_1A^+)+...+w_m(c^+-Q_mA^+)=w_1c^{\land}_1+...+w_mc^{\land}_m$.\\
	\end{center}

	Since $D$ is the same, (6.2) holds, \;  \textbf{End}
	
	It is obvious.\\
	\textbf{Corollary 6.1} If $F^*=w_1F_1+...+w_mF_m$, then
	\begin{equation}\label{6.3}
 (F^*)\downarrow=w_1(F_1)\downarrow+...+w_m(F_m)\downarrow.
	\end{equation}
	
	Set that:
	\begin{center}
	$[F_1,...,F_m]|b=\{y|\exists [F_1,...,F_m]\},$ 
	\end{center}

	Which is the column generated by$[F_1,...,F_m]$, it is not difficult to prove that.\\
	\textbf{Corollary 6.2} Let $P^*$ be the optimal point in $D$, if $P^*$ belongs to $[F_1,...,F_m]|b$, then
	there is a weight vector $w^*$ such that
	\begin{equation}\label{6.5}
 (F^*)\downarrow=(w^*_1F_1+...+w^*_mF_m)\downarrow=P^*.
	\end{equation}
	
	Thus, the solution of linear programming is transformed into an optimization problem of
	variable weight combination.
	
	We have no ability to solve the problem in the paper but presents a brief idea as follows:
	
	Set $r=(b,b)$, called the falling rate of parameter $t$. For example, if $b=(4,1,4,2,6)$, then $r=73$, It means that the $t$-value plus 1, the falling height will plus 73. Inversely, to make the falling height go down 1, the $t$-value should plus
	$1/r$. In the example, to fall down 2.6 it needs $t$ decrease $2.5/73=0.036$. Without the rate $r$, we can’t do good optimization on $t$-value tuning.\\
	\textbf{Example 6.2} Based on the 5 feasible vertices in Example 6.1, the variable weight combination is sought to make its landing point as close to the optimal point as possible.
		\begin{center}
			Step 1 Writing $t$-value vectors of $F_1-F_5$ and $\underline{F}$ as follows:
		\begin{tabular}{cccccccccccc}
			&$x_1$ &$x_2$ &$x_3$ &$x_4$ &$x_5$ &$x_6$ &$y_1$ &$y_2$ &$y_3$ &$y_4$ &$y_5$ \\ 
			$t_1$:&  0.08 & 4.3&&  0.2&  0.44&  0.16&  0.88 & 0.2 & 0&  0&  0\\ 
				$t_2$: &0 &9.5& &0.07& 0.28& 0.54& 0& 2.6 &  0 &3& 0\\ 
				$t_3$:&  0.05& 6.7&&  0.03&  1.03&  0&  0.48& 0.1 & 0.08&  1.43&  0\\ 
				$t_4$: &0.08 &6.6&& 0& 1.07& 0& 0.48& 0.1& 0&  1.72 & 0\\ 
			$t_5$:&  0.08 & 6.5&&  0.02&  0.99&  0.02&  0.47 & 0.03 & 0.05&  0.43&  0.03\\ 
				$\underline{t}$:& 0.06& 6.67& & 0.07 &0.9& 0.05&0.48 &0.6 & 0.03&1.5 &0.01\\ 
		\end{tabular} 
	\end{center}
	
In these 6 vectors, $\underline{t}$ is the best, since it does not include zero. Then ask where are its minimum and second minimum. The last one corresponding to $y_5$ is the minimum 0.01, the second minimum 0.018 corre sponds to $y_3$. Which
vector has biggest value corresponding to $y_5$? The vector $t_5$ has 0.05 there, unfortunately, which has zero
corresponding to $y_3$. There needs to do weight combination. 

If we take weights (0.8,0.2): then
0.01$\times$0.8+0.05$\times$0.2=0.018, the minimum is increasing but 0.018$\times$0.8+0$\times$0.2=0.0144, the second minimum becomes
minimum. So we can improve the weights as (0.862 0.138):\begin{center}
$0.01\times 0.862+0.005\times 0.138=0.0155$,\\
$0.018\times 0.862+0\times 0.138=0.0155$.
\end{center}

This is the better result. Then we get the weighted center of horizontal feasible outline:\begin{center}
$F^*=0.862 F_5 +0.138 \underline{F} =0.028 (F_1+F_2+F_3+F_4)+0.89 F_5=(1.8, 0.2, 0.06, 2.7, 0.26)$.\\
$F^*\downarrow=F^*-0.0155 b=(1.74, 0.18, 0, 2.67,0.17)=13.5.$\\
$h(F^*)-h(F^*\downarrow) =14.6-13.5=1.1$.
\end{center}

The effectove of variable weight  ombination falling is 1.1/2.6=42\%, which is
better than the feasible central falling. However, it is not a good result in the example. Indeed, if the optimal point is within the falling region of [$F_{\Delta}$], then, we can make the landing point directly attach the optimal solution. However, we can return to the Step 1, but we do not state the process in detail.
	\section{Tri-Skill combination}
	Tri-skill stands for 1. Highest principle on cutter-selection; 2. New column elimination; 3. Horizontal feasible central falling/sliding. These three algorithms all have their own advantages, but the strong power is in the combination of the three skills:
	
	Given a Simplex table \ref{Tab2.1}, Assume $b\geq0$; And given a feasible point $F_o$;

	Input index set of constraint planes $J:=\{1,...,n,n+1,...,n+m\};$\\
	Step 1: Check if the cutting degree $c^{\land}$ is non-positive. If so, stop and set the cone vertex as a dual optimal point.\\
	Step 2: Check whether slope vector $s=B^{-1}b$ has same terms, If so, make transformation to occur horizontal edges, and check feasible intervals on each horizontal edge. The convex set spanned by these feasible intervals is the set of optimal points; If there is no feasible points on all horizontal edges, they can be all cut off, and maintains the strict normal cone;\\
	Step 3: $F^{k+1}:=F^k\downarrow$ (doing on current table \ref{Tab2.2}), check if the relative error $\epsilon=(h(V)-h(F^{k+1}))/h(V)$ is enough small, you can tuning an error threshold value, $\epsilon<0.1$, for example, if it is not, calculate the supportor of horizontal cover, doing fesible central falling again and again; Else, go to Step 4; \\
	Step 4: Do pivoting according to the highest algorithm, recording critical heights in the bottom of negative cutting degree columns. Deleting the eliminatable planes’ indices from index set $J$. Repeat step 4 within $m$ times. After $m$ times, go back to step 1.
	
	The contribution of the algorithm: 1. In big data LP, optimal points always occur in horizontal edges, the algorithm provides a fast way to face the new challenge. 2. Doing feasible falling in the start, by a few steps (within 10 times), attaching the golden cone or near by that, then employ the highest cutting plus column elimination, the linear programming is not a complex but a ordinary calculation.
	\section{Conclusions} \label{sec:Conclusions}
	Linear programming is of extreme significance to artificial intelligence and data science. It is not only an indispensable computing tool, It is an important topic of concern to the theory of factor space\cite{15,16,17}. The realization of strong polynomial algorithm for linear programming has solved a major theoretical crisis for the development of artificial intelligence.
	\section*{\centerline{Acknowledgment}}
	Thanks to Liaoning Technical University, Big Data center of CAS for their supports; thanks to Professors and scholars Yixin Zhong, Huacan He, Zhende Huang, Yingjie Tian, Yongyi Chen, Kaiqi Zou, Hongxing Le, He Ouyang, Liyan Han, Xuehai Yuan, Qing He, Fanhui Zeng, Tiejun Cui, Yanke Bao, Haitao Liu, Jianjun Li, Jianwei Guo, Linjie Pu and Runjun Wan for their Supports and helps with enthusiasm.
	
\end{document}